\documentclass[11pt]{article}
\usepackage{graphicx}
\usepackage{amssymb,amsfonts,amsthm}
\usepackage{amsmath,amsthm,amssymb}
\usepackage{amsfonts}
\newcommand{\LDC}{{\mathrm{LDC}}}
\newcommand{\Con}{{\mathrm{Con}}}

\newcommand{\dist}{{\mathrm{dist}}}

\newtheorem{theorem}{Theorem}[section]
\newtheorem{lemma}[theorem]{Lemma}

\theoremstyle{definition}

\newcounter{ppp}

\newcommand{\la}{\langle}
\newcommand{\ra}{\rangle}

\newcommand{\area}{{\rm area}}

\newcommand{\vk}{van Kampen }

\begin{document}

\renewcommand{\theequation}{\thesection.\arabic{equation}}

\title{A finitely presented group with two non-homeomorphic asymptotic cones}

 \author{A.Yu. Ol'shanskii, M.V. Sapir\thanks{Both authors were supported in part by
the NSF grant DMS 0245600. In addition, the research of the first
author was supported in part by the Russian Fund for Basic Research
02-01-00170,  the research of the second author was supported in
part by the NSF grant DMS 9978802 and the US-Israeli BSF grant
1999298.}}
\date{}
\maketitle

\begin{abstract}
We give an example of a finitely presented group $G$ with two
non-$\pi_1$-equivalent asymptotic cones.
\end{abstract}

\section{Introduction}

Asymptotic cones of groups were introduced by Gromov \cite{Gr1} to
prove that a group of polynomial growth is virtually nilpotent. In
\cite{DW}, the concept was generalized to arbitrary finitely
generated groups. By definition, an asymptotic cone of a group
depends on the choice of an ultrafilter and a choice of an
increasing sequence of numbers $d_n\to \infty$. Nevertheless, in
many cases all asymptotic cones of a group turn out to be
homeomorphic. In particular, this is the case for hyperbolic groups
\cite{DP}, for nilpotent groups \cite{Gr1}, \cite{Pansu}, etc. In
\cite{Gr3} (see Question 2.$B_1$ (c)), Gromov asked whether there
exists a finitely generated (finitely presented) group with two
non-homeomorphic asymptotic cones.

S. Thomas and B. Velicovic \cite{TV} gave an example of a finitely
generated group $H$ with two non-homeomorphic asymptotic cones. C.
Drutu and M. Sapir \cite{DS} gave an example of a finitely generated
group with continuum pairwise non-homeomorphic (and even
non-$\pi_1$-equivalent) asymptotic cones. On the other hand, by L.
Kramer, S. Shelah, K. Tent and S. Thomas \cite{KSTT}, if the
Continuum Hypothesis is true then continuum is the maximal number of
non-isometric asymptotic cones a finitely generated group can have.
If the Continuum Hypothesis is not true, they give an example of a
finitely presented group with $2^{2^{\aleph_0}}$ pairwise
non-homeomorphic asymptotic cones (if the Continuum Hypothesis is
true, their group has unique asymptotic cone).

It is essential in the proofs in \cite{TV} and \cite{DS} that the
groups in both papers are limits of hyperbolic groups but
non-hyperbolic themselves (all asymptotic cones of non-elementary
hyperbolic groups are isometric \cite{Gr3}, \cite{DP}). Therefore
these groups are not finitely presented. Moreover, some of the
asymptotic cones of these groups are $\mathbb{R}$-trees, so they
cannot be finitely presented because it has been proven by M.
Kapovich and B. Kramer \cite{KK} that if an asymptotic cone of a
finitely presented group is an $\mathbb{R}$-tree then the group is
hyperbolic and so all its asymptotic cones are isometric.

The question of whether for a {\em finitely presented} group
asymptotic cones can be non-ho\-meo\-mor\-phic (non-isometric)
independently of the Continuum Hypothesis was open. The goal of this
note is to give a positive answer to this question.

We use Dehn functions of groups. In \cite{Ol}, the first author
constructed (using some ideas from \cite{OS}) a finitely presented
group $G$ whose Dehn function $f(n)$ satisfies the following two
properties:

(P1) there are sequences of positive numbers $d_i\to\infty$ and
$\lambda_i\to\infty$ such that $f(n)\le c n^2 $ for arbitrary
integer $n\in \cup_{i=1}^{\infty} [\frac{d_i}{\lambda_i} ,\; \lambda
_i d_i]$ and some constant $c$,

(P2) there is a positive constant $c'$ and an increasing sequence of
numbers $n_i\to \infty$ such that $f(n_i)/n_i^2\to \infty$ but for
every $i$, and for every integer $n$ with $n\le c' n_i$, we have
$f(n)\le cn_i^2$.

\begin{theorem} \label{th2} Let $G$ be a finitely presented group
satisfying (P1) and (P2). Then $G$ has two asymptotic cones, one of
which is simply connected and another one is not.
\end{theorem}

It is known \cite{Pap} that if the Dehn function of a group is
quadratic then all its asymptotic cones are simply connected. We
slightly modify Papasoglu's argument and show that property (P1)
implies that one of the asymptotic cones of $G$ is simply connected.

On the other hand, by a result of Gromov \cite{Gr3} (see also
\cite{Dr}), if the Dehn function of a group is not bounded by a
polynomial then the group has a non-simply connected asymptotic
cone. The Dehn function of $G$ is bounded by $n^3$ (in fact, by
$n^2\log n/\log\log n$) but Property (P2) allows us to apply
essentially Gromov's argument and show that $G$ has a non-simply
connected asymptotic cone.

Note that the group $G$ from \cite{Ol} is an $S$-machine in the
terminology of the second author, so it is a multiple HNN extension
of a free group with finitely generated associated subgroups (see
\cite{SBR}, \cite{OS}). In particular, $G$ has cohomological
dimension 2.

\section{Proof}

Recall the definition of an asymptotic cone.  A non-principal
ultrafilter $\omega$ is a finitely additive measure defined on all
subsets $S$ of ${\mathbb N}$, such that $\omega(S) \in \{0,1\}$ and
$\omega(S)=0$ if $S$ is a finite subset. For a bounded function $f:
{\mathbb N}\to {\mathbb R}$ the limit $\lim_{\omega} f(i)$ with
respect to $\omega$ is the unique real number $a$ such that
$\omega(\{i\in {\mathbb N}: |f(i)-a|<\epsilon\})=1$ for every
$\epsilon>0$.

Let $(X,\dist)$ be a metric space. Fix an arbitrary $x_0\in X$, and
a sequence of {\em scaling constants} $d_i\to \infty$.  Consider the
set of sequences $g\colon {\mathbb N}\to X$ such that
$\dist(f(i),x_0) \le c d_i$ for some constant $c=c(f)$. Two
sequences of this set ${\cal F}$ are said to be equivalent if
$\lim_{\omega}\frac{\dist(f(i),g(i))}{d_i} =0$. The asymptotic cone
$\Con^{\omega}(X, (d_i))$ is the quotient space ${\cal F}/\sim$
where the distance between equivalence classes $[f]$ and $[g]$ is
equal to $\lim_{\omega}\frac{\dist(f(i),g(i))}{d_i}$. The asymptotic
cone is a complete space; it is a geodesic metric space if $X$ is a
geodesic metric space (\cite{Gr3}; \cite{Pap}). If $[f]$ is an
element of $\Con^\omega(X,(d_i))$ then we say that $f(i)$ {\em
converges} to $[f]$. Note that $\Con^\omega(X, (d_n))$ does not
depend on the choice of $x_0$.

An {\em asymptotic cone} of a group $G$ with a word metric is
isometric to the asymptotic cone of its Cayley graph (considered as
the 1-skeleton of the Cayley complex). The asymptotic cones of the
same group relative to two finite generating sets are bi-Lipschitz
equivalent. Note that in \cite{Gr3}, \cite{Pap} and other papers, a
more restrictive definition of asymptotic cone was used: it was
always assumed that $(d_i)=(i)$. It was observed in \cite{Ri},
however, that if, say, all $d_i$'s are different integers then
$\Con^\omega(G,(d_i))$ is isometric to a cone
$\Con^{\omega'}(G,(i))$ for some $\omega'$. Since in all the
asymptotic cones considered in this paper, all $d_i$'s are different
integers, they are isometric to restricted asymptotic cones.

As in \cite[p. 792]{Pap}, we define an $n$-{\em gone} $P$ in a
geodesic metric space $(X, \dist)$ as a map from the set of vertices
of the standard regular $n$-gone $\bar S_n$ in the plane into $X$.
If $X$ is a Cayley graph of a group, we shall always assume that
elements of $P$ are vertices of the graph, i.e. they belong to $G$.
An {\em side} of $P$ is a pair of vertices of $P$ corresponding to
the pair of vertices connected by an edge in $\bar S_n$; the {\em
length} of a side is the distance between these two vertices. The
{\em perimeter} (length) of $P$ is the sum of lengths of its sides.
A {\em partition} of $\bar S_n$ is  a collection of discs
$D_1,...,D_k$ such that $\bar S_n=\partial(D_1\cup...\cup D_k)$ and
$D_i\cap D_i\subseteq \partial D_i\cap \partial D_j$ if $i\ne j$. A
{\em vertex} of a partition $D_1,...,D_k$ is either a vertex of
$\bar S_n$ or a point on $\partial D_1\cup...\cup \partial D_k$ such
that for every open set $U$ containing this point, the intersection
$U\cap \bigcup \partial D_i$ is not homeomorphic to an interval.

A {\em partition} of $P$ is a map $\Pi$ from the set of vertices of
a partition of $\bar S_n$ into $X$ taking the vertices of $\bar S_n$
to $P$. {\em Vertices} of the partition of $P$ are images of
vertices of the partition of $\bar S_n$ under $\Pi$. Note that for
each $i=1,...,k$, the images of the set of vertices of belonging to
$\partial D_i$ form a polygon in $X$. This polygon will be called a
{\em piece} of the partition $\Pi$.

\begin{lemma} \label{793} Let $(X, \dist)$ be a geodesic metric space,
and $P$ a polygon in $\Con^{\omega}(X, (d_i))$ with vertices
$P_1,\dots,P_n$. Assume that $P$ satisfied the following Loop
Division Condition:

$\LDC(k)$: There exists a sequence of polygons
$Q^i=(P_1^i,\dots,P_n^i)$ in $X$ such that $P_j^i$ converges to
$P_j$ for $j=1,\dots,n$, and every $Q^i$ can be partitioned into $k$
pieces whose perimeters are less than or equal to $\frac12
perimeter(Q^i)$.

Then the polygon $P$ can be partitioned in $\Con^{\omega}(X, (d_i))$
into $k$ pieces whose perimeters do not exceed $\frac12$ of the
perimeter of $P$.

\proof This assertion is proved in \cite{Pap}. (See the proof of the
Proposition formulated on page 793; that proof works without
changes, although the formulation of Lemma \ref{793} slightly
differs from the formulation of the cited Proposition.) \endproof
\end{lemma}

\begin{lemma} \label{793b} Let $X$ be a complete geodesic metric space such
that for some integer $k$, every polygon $P$ of the asymptotic cone
$\Con^{\omega}(X, (d_i))$ satisfies $\LDC(k)$. Then $X$ is simply
connected.
\end{lemma}
 \proof Again, it suffices to repeat the proof of the Proposition formulated
on the bottom of page 793 of \cite{Pap} (though our formulation
differs from that in \cite{Pap}, and one should refer to Lemma
\ref{793} now). \endproof

The following version of Papasoglu's lemma is now formulated for
arbitrary planar triangular map, i.e. for a map whose faces are of
(combinatorial) perimeter at most 3.

\begin{lemma}\label{Pa} Let $\Delta$ be a  triangular map whose perimeter
 $n$ is at least 200. Assume
that the area of $\Delta$ does not exceed $Mn^2$. Then there is $k$
depending on $M$ only, such that
$\Delta=\Gamma_1\cup\dots\cup\Gamma_k$ where $\Gamma_i$, ($i=1,\dots
k$) are submaps of $\Delta$, and $\Gamma_i\cap \Gamma_j$ ($0\le
i<j\le k$) is empty or a vertex, or a simple path, and perimeter
$|\partial\Gamma_i|$ is at most $n/2$ for all $i=1,\dots,k$.
\end{lemma}

\proof The proof of the Theorem formulated in \cite{Pap}, page 799,
does not use the labels of diagram edges, and so it also works for
maps. Although it is assumed in \cite{Pap}, that $\area (\Delta)\le
M|\partial\Delta|^2$ for {\em all} minimal van Kampen diagrams over
a triangular group presentation, the proof uses this quadratic
isoperimetric inequality only for one diagram $\Delta$.  The
assertion of Lemma \ref{Pa} is therefore correct.
\endproof

Let $Q$ be a polygon in the Cayley graph of $G$. Connect the
vertices of each side of $Q$ by a geodesic, then the product of
labels of these geodesics viewed as a cyclic word is called a {\em
label} of $Q$ (a label depends on the choices of the geodesics, of
course).

\begin{lemma}\label{simcon} Let $f$ be the Dehn function of a finite
group presentation $G = \la A \mid R\ra$ satisfying (P1). Then the
asymptotic cone $\Con^{\omega}(G, (d_i))$ is simply connected for
arbitrary non-principal ultrafilter $\omega$.
\end{lemma}

\proof Let $P=(P_1, \dots, P_m)$ be a polygon in
$\Con^{\omega}(X,(d_i))$ with pairwise distinct vertices $P_1,\dots,
P_m$. Consider a sequence of polygons $Q^i=(Q_1^i,\dots,Q_m^i)$ in
the Cayley graph $\Gamma(G,A)$ of the group $G$, such that $Q_j^i$
converges to $P_j$ for every $j$. Since $\dist(P_j,P_{j'})>0$ in the
cone for $j\ne j'$, there are constants $\alpha$ and $\beta$
independent of $i$ and $j$ such that $l_i=\mathrm{perimeter}(Q_i)\in
[\alpha d_i, \beta d_i]$ for almost all $i$-s (with respect to
$\omega$). Property (P1) implies that
\begin{equation}\label{vnutri}
\omega (I)=1 \;\; for\;\; I=\{i\mid l_i\in [\frac{d_i}{\lambda_i}
,\; \lambda_i d_i]\}
\end{equation}

Van Kampen's Lemma provides us with a minimal diagram $\Delta_i$
over the presentation $G=\la A\mid R\ra$ such that the boundary
label of $\Delta_i$ is a label of the polygon $Q^i$ in the Cayley
graph $\Gamma(G,A)$, $i\in I$.  By Lemma \ref{Pa} and formula
(\ref{vnutri}), there is a constant $k=k(c)$ such that $\Delta_i$
can be partitioned into subdiagrams $\Gamma_1^i,\dots,\Gamma_k^i$
with perimeters at most $l_i/2$. Then the polygon $Q^i$ can be
accordingly partitioned into discs $D_1^i,\dots,D_k^i$ in the Cayley
graph, and the perimeters of these discs do not exceed $l_i/2$.
Hence, in the cone, every polygon $P$ satisfies $\LDC(k)$, and, by
Lemma \ref{793b}, the cone $\Con^{\omega}(G, (d_i))$ is simply
connected.
\endproof

\begin{lemma}\label{dr}
Let $f$ be the Dehn function of a finite group presentation $G=\la
A\mid R\ra$ satisfying property (P2). Then, for arbitrary
non-principal ultrafilter $\omega$, the asymptotic cone
$\Con^{\omega}(G, (n_i))$ is not simply connected.
\end{lemma}

\proof Property (P2) implies existence of a positive constant $c'<1$
such that $f(n_i)/f(c'n_i)\to \infty$.

Assume there is a number $k$ such that for every $i$, an arbitrary
polygon $Q$ of the Cayley graph $\Gamma (G,A)$ of length $l$,
$c'n_i\le l\le n_i$, can be partitioned into at most $k$ pieces of
perimeter $l/2$. It follows that every loop of length $n_i$ can be
partitioned into at most $K=k^{1-\log_2 c'}$ pieces of perimeter at
most $c'n_i$.

Let $\Delta_i$ be a van Kampen diagram with perimeter $n_i$ and area
$f(n_i)$ that has minimal area among all diagrams with the same
boundary label. For each $i$, consider the loop $Q^i$ in the Cayley
graph $(\Gamma, \dist)$ of $G$, whose label is equal to the boundary
label of $\Delta_i$'s. It follows from our assumption that for every
$i$, there is a partition of $Q^i$ into at most $K$ pieces  of
perimeter $\le c'n_i$. The smallest area \vk diagram having the same
label as $Q^i$ has area at most $f(c'n_i)$. Therefore the area of
the minimal diagram $\Delta$ cannot exceed $Kf(c'n_i)$. Hence
$f(n_i)/f(c'n_i)\le K$ for every $i$. This contradicts property
(P2).

Therefore our assumption was false, and there is no such number $k$.
Also there is a constant $c_0$ such that the radius of $Q^i$ (i.e.
$\max (\dist(x,y)\mid x,y\in Q^i$) is at least $c_0n_i$. Otherwise
we could easily partition $Q^i$ into a bounded number of loops with
length $\le c'n_i$ (which can be ruled out as in the previous
paragraph). Consider the $\omega$-limit of the sequence of (finite)
sets $Q^i$, i.e. the set of all elements $[(x_i)]$ where $x_i\in
Q^i$. It is easy to see (cf., for example,\cite{Buragos}) that the
$\omega$-limit of $Q^i$ is a loop $P$ in $\Con^\omega(G, (n_i))$ of
length at least $c_0$. Indeed, one can parametrise each loop $Q^i$
by its arc length by a function $x_i\colon [0,1]\to (\Gamma,
\dist/{n_i})$, then $P$ has parametrization $x\colon [0,1]\to
\Con^\omega(G,(n_i))$ where $x(t)=[(x_i(t))]$ for each $t\in [0,1]$.

The loop $P$ has no finite partition into pieces $P_1,P_2,\dots,
P_k$ whose perimeters do not exceed a half of the perimeter of $P$.
Therefore the loop $P$ is not contractible. (For more details
justifying the last two phrases, see \cite{Gr3}, \cite{Pap} or the
proof of Theorem 4.4 in \cite{Dr}.)
\endproof

{\bf Proof of Theorem \ref{th2}.} The theorem follows from lemmas
  \ref{simcon}, \ref{dr}.

\medskip

{\bf Acknowledgement.} The authors are grateful to Cornelia Drutu,
Denis Osin and Panos Papasoglu for helpful conversations.

\begin{minipage}[t]{3 in}
\noindent Alexander Yu. Ol'shanskii\\ Department of Mathematics\\
Vanderbilt University \\ alexander.olshanskiy@vanderbilt.edu\\
 and\\ Department of
Higher Algebra, MEHMAT\\
 Moscow State University\\
olshan@shabol.math.msu.su\\
\end{minipage}
\begin{minipage}[t]{3 in}
\noindent Mark V. Sapir\\ Department of Mathematics\\
Vanderbilt University\\
m.sapir@vanderbilt.edu\\
http://www.math.vanderbilt.edu/$\sim$msapir\\
\end{minipage}

\end{document}